\documentclass[oneside,english]{amsart}
\usepackage[latin1]{inputenc}
\usepackage{amssymb}

\makeatletter

\title{Quasi-factors for infinite-measure preserving transformations}
\author{Tom Meyerovitch}
\newtheorem{theorem}{Theorem}[section]
\newtheorem{lemma}[theorem]{Lemma}
\newtheorem{proposition}[theorem]{Proposition}
\newtheorem{corollary}[theorem]{Corollary}

\input xy
\xyoption{all}

\usepackage{color}
\usepackage{epsfig}

\newcommand{\hide}[1]{}
\usepackage{babel}
\makeatother

\newcommand{\X}{\mathbb{X}}

\newcommand{\B}{\mathcal{B}}
\newcommand{\hrel}{{h_f}}
\newcommand{\cmpts}{conservative, measure preserving system}

\begin{document}

\begin{abstract}
This paper is a study of Glasner's definition of quasi-factors in
the setting of infinite-measure preserving system.
The existence of a system with zero Krengel entropy and a
quasi-factor with positive entropy is obtained. On the other hand,
relative zero-entropy for conservative systems implies relative
zero-entropy of any quasi-factor with respect to its natural
projection onto the factor. This extends (and is based upon) results
of Glasner, Thouvenot and Weiss
\cite{glasner_thouvenot_weiss_2000,glasner_weiss95}. Following and
extending Glasner and Weiss \cite{glanser_wiess_02}, we also prove
that any conservative measure preserving system with positive
entropy in the sense of Danilenko and Rudolph
\cite{danilenko_rudolph07} admits any probability preserving system
with positive entropy as a factor. Some applications and connections
with Poisson-suspensions are presented.
\end{abstract}

\email{tomm@post.tau.ac.il}


\maketitle

\section{Introduction}
Glasner \cite{glasner_quasifactors_1983} introduced the definition of a \emph{quasi-factor}
of a probability preserving dynamical system, as a conceptual generalization of a factor system.
A quasi-factor is a probability-preserving action on the space of measures, which satisfies a barycenter condition.
As proved in \cite{glasner_quasifactors_1983}, any
factor of a probability preserving system gives raise to a quasi-factor via the disintegration of $\mu$
over $\nu$. In Glasner's work and follow-ups with Weiss
it has been shown that certain properties of a dynamical systems are reflected by similar properties of their quasi-factors:
For example, any quasi-factor of a Kronecker system is again a Kronecker system,
and any quasi-factor of a zero-entropy system has zero entropy. Contrary to this, it has been proved
by Glanser and Weiss in \cite{glanser_wiess_02}
that a system with positive entropy admits any other positive entropy system as a quasi-factor.


In the present paper we extend the notion of quasi-factors for $\sigma$-finite measure preserving systems.
A main difference in the theory of quasi-factors for infinite-measure systems is the following:
A quasi-factor of an infinite measure preserving system is a probability preserving system.
Thus, when studying quasi-factors we attempt to learn about infinite measure preserving systems via
some associated probability preserving systems.

Our motivation for this comes partly from the role of Poisson-suspensions in infinite ergodic theory.
Research of Poisson-suspensions was carried out by Roy \cite{roy_thesis,roy_poisson}. Questions regarding entropy of Poisson suspensions were addressed  in \cite{jmrr_poisson_07}.
We observe that Poisson-suspensions are naturally interpreted as quasi-factors.
As a consequence of theorem \ref{thm:quasi-factor-relative-entropy} below,
 relative zero-entropy for conservative systems implies relative zero-entropy of Poisson
 suspensions,
a result obtained by direct methods in \cite{jmrr_poisson_07}.

\textbf{Acknowledgments:} I like to thank Eli Glasner and Benji
Weiss for interesting remarks and discussions. This work is a part
of my Ph.D in Tel-Aviv University. I owe special thanks to Jon
Aaronson, my Ph.D advisor for his valuable guidance,  and infinite
patience. I gratefully acknowledge the support of support of the
Crown Family Foundation Doctoral Fellowships ,USA.

\section{ Preliminaries and definitions}
Let $(X, \mathcal{B},\mu)$ be a standard measure space.
Denote by $X^*$ the set of measures on $(X,\mathcal{B})$, and $\mathcal{B}^*$ the $\sigma$-algebra defined by: \begin{equation}
\label{eq:B_star_sigma_algebra}
\mathcal{B}^* = \sigma\left(\left\{ \{\gamma \in X^*:~ \gamma(B) \in [a,b]\}: B \in \mathcal{B},\, 0\le a \le b \le \infty ~\right\}\right)
\end{equation}

Any $\B$-measurable map $T:X \to X$  gives raise to a $\B^*$-measurable map $T_*:X^* \to
X^*$ by $T_*\gamma = \gamma \circ T^{-1}$. A \emph{quasi-factor} of $(X, \mathcal{B},\mu)$
 is a probability preserving dynamical system of the form $(X^*,\B^*,\xi,T_*)$ where
  $\xi$  is  a $T_*$-invariant probability on $X^*$, and whose barycenter is $\mu$:

\begin{equation}\label{eq:barycenter} \mu = \int_{X^*}\gamma d\xi(\gamma)
\end{equation}
By this we mean that for any $f \in L^1(X,\mathcal{B},\mu)$, the
expectancy under $\mu$ of the (random) integral of $f$ is the
integral of $f$ with respect to $\mu$. In other words, $\int f(x)
d\mu(x) = \int f^*(\gamma) d\xi(\gamma)$, with $f^*X^* \to
\mathbb{R}$ denoting the function defined by $f^*(\gamma)=\int
f(x)d\gamma(x)$. For equation \eqref{eq:barycenter} to make sense,
we require in particular that for each set $A \in \mathcal{B}$, if
$\mu(A)=0$ then the set $\{\gamma \in \X^* :~ \gamma(A)>0\}$ is null
modulo $\xi$.

It may not immediately be evident to the reader that a probability space $(X^*,\mathcal{B}^*,\xi)$ as described above is a standard probability space,
in the sense that it is isomorphic to a completion of a Borel-probability measure on a Polish space.
 We will immediately make this clear, but first,
let us recall Glasner's original definition of a quasi factor from \cite{glasner_quasifactors_1983}:

Assume that $(X,\mathcal{B},\mu,T)$  is a probability preserving transformation, compatible with a certain topological structure.
$X$ is a compact metric space, $\mathcal{B}$ is the completion of the Borel $\sigma$-algebra and $T:X \to X$ is a homeomorphism of $X$.
Glasner defined a quasi-factor of $(X,\mathcal{B},\mu,T)$ to be any probability preserving system of the form $(X^\dagger,\mathcal{B}^\dagger,\xi,T_*)$
where $X^\dagger$ denotes the set of Borel probability measures on $X$,
equipped with the compact weak-$*$ topology, $\mathcal{B}^\dagger$ is the Borel-$\sigma$ algebra of $X^\dagger$,
 $T_*:X^\dagger \to X^\dagger$ is the homeomorphism of $X^\dagger$
defined by $T_* \nu =\nu \circ T^{-1}$ and $\xi$ is a probability measure on $X^\dagger$ whose barycenter is $\mu$. This means that for every continuous function $f:X \to \mathbb{R}$ we have the inequality
$\int f(x)d\mu(x) = \int f^*(\gamma) d\xi(\gamma)$.

Glanser has proved that his definition is independent of the topological structure.

The fact that any quasi-factor $(X^*,\mathcal{B}^*,\xi)$ is a standard probability space follows from this lemma:
\begin{lemma}\label{lem:standard_space_quasifactor}
  Let $(X,\mathcal{B},\mu)$ be a standard measure space,
  and let $(X^*,\mathcal{B}^*,\xi)$ be a probability space satisfying the barycenter condition~\eqref{eq:barycenter}.
  Let $\alpha=\{A_n\}_{n=1}^{\infty}$ be a measurable partition of $X$ into sets with $\mu(A)<\infty$,
  then $(X^*,\mathcal{B}^*)$ is isomorphic $\mod \xi$ to
  $$\left((\prod_{n=1}^{\infty} C A_n^\dagger, \bigotimes_{n=1}^{\infty}\mathcal{B}(C A_n^\dagger)\right),$$
  where $C X:= (X \times \mathbb{R}_+) / (X \times \{0\})$ denotes the cone of a space $X$.

  Here $A_n^\dagger$ denotes the set of probabilities on $A_n$, and $\mathcal{B}(A_n)^\dagger$ is the natural $\sigma$-algebra  on $A_n^\dagger$
  which comes from $\mathcal{B} \cap A_n$.
\end{lemma}
\begin{proof}
  The barycenter condition grantees that $\xi$-almost-surely the measure of every $A_n$ is finite.
  Consider the map $\Phi:X^* \to \prod_{n=1}^{\infty} C A_n^\dagger$ defined by
  $$\Phi(\gamma)_n=\left\{ \begin{array}{cc}
  \left(\frac{\gamma \mid_{A_n}}{\gamma(A_n)},\gamma(A_n)\right) & \mbox{if } \gamma(A_n)>0\\
  0 & \mbox{if } \gamma(A_n)=0
  \end{array} \right.$$
  The map $\Phi$ is well defined on a set of full $\xi$-measure, it is $\mathcal{B}^*$-measurable and injective $\mod \xi$.
\end{proof}
As $(A_n^\dagger,\mathcal{B}(A_n)^\dagger)$ is clearly a standard measurable space, so is $(X^*,\mathcal{B}^*,\xi)$.

The consistency  of  our definition with  Glasner's is now clear: If $\mu$ is a probability measure and $\xi$ is an ergodic
quasi-factor of $\mu$, then $\gamma(X)$ is constant $\xi$-a.e., and
therefore $\gamma(X)=1$ $\nu$-a.e. Thus, for ergodic quasi-factors,
our definition of a quasi-factor coincides with Glasner's
original definition from \cite{glasner_quasifactors_1983} for the probability
preserving category.


\subsection{Measure preserving and  nonsingular-factors of measure preserving systems}

A measure preserving system $(Y,\mathcal{C},\nu,S)$ is a \emph{factor} (or homomorphic image) of $(X,\mathcal{B},\mu,T)$
if there is a measurable $\pi:X \to Y$ which sends $\mu$ to $\nu$ and intertwines $T$
with $S$: $S\circ \pi = \pi \circ T$ and $\nu = \mu \circ \pi^{-1}$. In such case we say that $T$ is an \emph{extension} of $S$, and $\pi$ is a
\emph{factor-map}.

Any quasi-factor $\xi$ of $(X,\mathcal{B},\mu,T)$
gives raise to a quasi-factor $\pi_*(\xi)$ of
a factor $(Y,\mathcal{C},\nu,S)$, as $\pi_*:X^* \to Y^*$ maps a measure on
$(X,\mathcal{B})$ to a measure on $(Y,\mathcal{C})$ and intertwines $T_*$ and $S_*$.

More generally, $(Y,\mathcal{C},\nu,S)$ is a \emph{non-singular} factor of $(X,\mathcal{B},\mu,T)$
if there is $\pi:X \to Y$ which intertwines $T$ and $S$ and such that $\nu$ has the same null-sets as $\mu$.
When $T$ and $S$ are both probability preserving and ergodic, these definitions coincide.

When $T$ and $S$ are measure-preserving, conservative and ergodic, $T$ is a non-singular factor of $S$
iff there exists a constant $c \in (0,\infty]$ such that $(Y,\mathcal{C},c\cdot\nu,S)$ is a measure preserving factor of $(X,\mathcal{B},\mu,T)$.
In this case, following Aaronson \cite{aaro_inf_book}, we say that $S$ is a $c$-factor of $T$. In this situation, if $\nu(Y)<\infty$ and $\mu(X)=\infty$,
we must have $c=\infty$, in which case $S$ is an $\infty$-factor of $T$.

\begin{proposition}\label{prop:factor_is_quasifactor}
  Let $(\Omega,\mathcal{F},p,\tau)$ be a probability preserving system which is a non-singular factor of a measure preserving system
  $(X,\mathcal{B},\mu,T)$. Then $(\Omega,\mathcal{F},p,\tau)$ is isomorphic to a quasi-factor of $(X,\mathcal{B},\mu,T)$.
\end{proposition}
\begin{proof}

  Let $\pi:X \to \Omega$  be a nonsingular-factor map.
There is a positive linear operator $\pi^\#:L^1(X,\mu) \to L^1(\Omega,p)$ of norm $1$ defined by the equation
$$ \int_A \pi^\#f(\omega)dp(\omega)= \int_{\pi^{-1}A} f(x) d\mu(x) ~~\forall A \in \mathcal{F}, \, f \in L^1(X,\mu)$$

This gives raise to a  disintegration of $\mu$ over $p$ by:
$$ \mu = \int_{\Omega}\mu_{\omega} dp(\omega)$$
where the conditional measures $\mu_\omega$ are defined by $\mu_\omega(B)= \pi^{\#}(1_B)(\omega)$ for $b \in \mathcal{B}$.

The map $\phi:\Omega \to X^*$ given by $\phi(\omega) = \mu_\omega$ is defined everywhere $\mod p$, is $\mathcal{F}$-measurable and injective.
Thus, $(X^*,\mathcal{B}^*,\phi_*(p),T_*)$ is a quasi-factor of $T$, and $\phi$ is an isomorphism of this quasi-factor to $(\Omega,\mathcal{F},p,\tau)$.
\end{proof}

Proposition \ref{prop:factor_is_quasifactor} appears for appears in
\cite{glasner_quasifactors_1983}, where the discussion was
restricted to the probability preserving context. This is the reason
for the name ``quasi-factor''.

\subsection{Entropy for $\sigma$-finite measure preserving transformation}
A most successful numeric invariant for  probability preserving
systems is the Kolmogorov-Sinai entropy. There have been several
works attempting to meaningfully extend this invariant for measure
preserving (or  non-singular) systems. Two definitions of
``entropy'' for infinite measure preserving transformations have
been suggested by Krengel and Parry: Krengel defines entropy of a
conservative measure preserving transformation is transformation
$(X,\mathcal{B},\mu,T)$ as follows \cite{krengel67}:
$$h_{Kr}(X,\mathcal{B},\mu,T)=\sup_{A \in
\mathcal{F}_+}\mu(A)h(A,\mathcal{B}\cap A,\bar{\mu}_A,T_A)$$ Where
$\mathcal{F}_+$ is the collection of sets in $\mathcal{B}$ with
finite, positive measure, $\bar{\mu}_A$ is the normalized
probability measure on $A$ obtained by restricting $\mu$ to
$\mathcal{B} \cap A$, and $T_A:A \to A$ is the induced map of $T$ on
$A$ defined by:
$$T_A(x) := T^{\phi_A(x)}(x)$$
with $\phi_A(x) := \min\{k \ge 1:~ T^k(x) \in A\}$ the \emph{first
return map} of $A$.

Parry \cite{parry_generators69} defines the entropy of a measure
preserving transformation by:
$$h(X,\mathcal{B},\mu,T)= \sup_{T^{-1}\mathcal{C} \subset
\mathcal{C}}H(\mathcal{C} \mid T^{-1}\mathcal{C}),$$

where the supremum  is over all $\sigma$-finite sub-$\sigma$-algebras
of $\mathcal{B}$ such that $T^{-1}\mathcal{C} \subset \mathcal{C}$.
 For probability preserving transformations, this definition coincides with standard definitions of Kolmogorov-entropy.

It is known that Krengel's entropy dominates Parry's entropy for any conservative transformation (see \cite{parry_generators69}).
It is unknown if these definitions coincide in general.

Given a factor map $\pi:X \to Y$ of systems $(X,\mathcal{B},\mu,T)$
and $(Y,\mathcal{C},\nu,S)$ the  \emph{relative entropy} $h(T \mid
S)$ can be defined in a natural way, extending the usual notion of
relative entropy of probability preserving systems. It has been
shown that this notion behaves well both with the definition of
Krengel. Danilenko and Rudolph studied relative entropy of
infinite-measure preserving systems  \cite{danilenko_rudolph07}.
They came up with following definition of entropy for a measure
preserving system:
$$\hrel(T) = \sup_{S}h(T \mid S),$$
where the supremum is taken over all factors of $T$.

\subsection{Relative Bernoulli extensions for infinite transformations}
Given a system $(Y,\mathcal{C},\nu,S)$ and two extensions $\pi_i:(X_i,\mathcal{B}_i,\mu_i,T_i) \to (Y,\mathcal{C},\nu,S)$, $i=1,2$,
we say $X_1$ and $X_2$ are \emph{relatively isomorphic} (over $Y$) if there exists an isomorphism $\phi:X_1 \leftrightarrow X_2$ such that
$\pi_1 = \pi_2 \circ \phi$, and such $\phi$ called is a relative isomorphism.

We begin with the following simple proposition, appearing in
\cite{aaron_keane_94} with a slightly different formulation:
\begin{proposition}\label{prop:rel_isomorphism_induced}
  Two extensions $\pi_i:(X_i,\mathcal{B}_i,\mu_i,T_i) \to (Y,\mathcal{C},\nu,S)$, $i=1,2$  of an ergodic \cmpts are relatively isomorphic iff
for some (any) $A \in \mathcal{C}$, the induced transformations of $T_1$ and $T_2$ on $\pi_i^{-1}A$ are relatively
isomorphic over the induced transformation $T_A$.
\end{proposition}
\begin{proof}
Obviously, if $\phi: X_1 \leftrightarrow X_2$ is a relative isomorphism over $S$ then
$\phi \mid \pi^{-1}_1A: \pi^{-1}_1 A \leftrightarrow \pi^{-1}_2A$ is a relative isomorphism of the induced transformations.

Suppose  $\phi_0:\pi^{-1}_1 A \leftrightarrow \pi^{-1}_2A$ is a relative isomorphism of the induced transformations.
Let $A_n = T^{-n}A \setminus \bigcap_{k=0}^{n-1}T^{-k}A$, then since $S$ is conservative and ergodic, $Y = \biguplus_{n=0}^{\infty}A_n$.
Define $\phi_n: \pi^{-1}_1 A_n \leftrightarrow \pi^{-1}_2A_n$ by $\phi_n := \phi \circ T^{n}$. It can now evident that the map $\phi:X_1 \to X_2$
defined by $\phi(x) = \phi_n(x)$ for $x \in \pi^{-1}_1 A_n$ is a relative isomorphism of $X_1$ and $X_2$.
\end{proof}

For probability preserving transformations, we say that $T$ is \emph{relatively Bernoulli} over $S$ if $(X,\mathcal{B},\mu,T)$ is isomorphic to
a system of the form $(Y\times \Omega, \mathcal{B}\otimes \mathcal{X},\nu \times p,T \times \sigma)$ where $(\Omega,\mathcal{X},p,\sigma)$
is a Bernoulli shift,
and the isomorphism $\phi: Y\times \Omega \to X$ respects $\pi$, meaning $\pi\circ\phi(y,z)=y$.

In the spirit of Krengel's definition of entropy, we extend this
notion: an ergodic \cmpts~ $T$ is relatively bernoulli over $S$ if
for some $A \in \mathcal{C}$ of finite measure, the induced
transforation $T_{\pi^{-1}A}$ is relatively Bernoulli over $S_A$.

To make sure the definition makes sense and extends the usual concept of relative Bernoulli, we prove the following lemma, making use of Thouvenot's
 relative isomorphism theory:
\begin{lemma}\label{lem:rel_bernoulli_induce}
With $T$ and $S$ as above,
let $A \in \mathcal{C}$ with $0< \nu(A)<\infty$.
$T$ is relatively Bernoulli over $S$ iff $T_{\pi^{-1}A}$ is relatively Bernoulli over $S_A$.
\end{lemma}

\begin{proof}
  Suppose  $T_A=S_A \times \sigma$, where $(\Omega,\mathcal{X},p,\sigma)$ is a Bernoulli shift with a process taking values in the countable set $\Sigma$,
  $\Omega= \Sigma^{\mathbb{Z}}$ and .
  Then the function $f: X \to \Sigma \cap \{0\}$ defined by
  $$f(y,\omega)= \left\{\begin{array}{cc}
  \omega_0 & y \in A\\
  0 & \mbox{otherwise}
  \end{array}\right.$$
  Determines a process on $T$ which is a relative generator with respect to $\mathcal{C}$.
  By Thouvenot's relative isomorphism theorem of \cite{thouvenot_relative_74} (see also \cite{kieffer_84}),
  since this process is $S$-relatively very weak-Bernoulli, it is relatively Bernoulli over $S$.

  For the other direction, suppose $T= S \times \sigma$ and $A \in \mathcal{C}$ is as above,
  let $f_A: A \times \Omega \to \Sigma^*$ be defined by $f_A(y,w)=(\omega_0,\ldots,\omega_{\phi_A(y)-1})$. Then the iterates of $f_A$ under $T_A$
   are jointly independent given $A \times \mathcal{X}$, and so $T_A$ is relatively Bernoulli over $S_A$.
\end{proof}

A naive way to define a Bernoulli extension of an infinite measure preserving system $S$ is being isomorphic to $S \times \Omega$.
The problem with this simple definition
is explained by the following proposition:
\begin{proposition}
If $B_1$ and $B_2$ are Bernoulli systems (not necessarily of the same entropy) and $T$ is a conservative, ergodic measure preserving transformation,
then   $X \times B_1 \cong X \times B_2$.
\end{proposition}
\begin{proof}
  By inducing onto a set of finite measure $A \in \mathcal{B}$, we see that both induced transformation are relatively
  Bernoulli over $T_A$ with infinite relative entropy, thus isomorphic. Applying proposition \ref{prop:rel_isomorphism_induced},
  we get an isomorphism for the original systems.
\end{proof}

\section{Some properties of quasi-factors}
Once we have introduced the definition of quasi-factors for infinite measure preserving transformations,
a basic question we can ask is the following:
Are there any properties of a measure preserving system which always pass to quasi-factors?

In the probability preserving category, Glasner has shown that any quasi-factor of a Kronecker system is isomorphic to a factor,
hence is also a Kronecker system.
Contrary to this, weakly mixing systems may have non weakly mixing systems as quasi-factors (see \cite{glasner_quasifactors_1983}).

In view of the results of \cite{glasner_weiss95} and \cite{glanser_wiess_02} once gets the impression that the
main property preserved by quasi-factors of probability preserving systems is having zero entropy. As we show in section \ref{sec:zero_entropy},
the situation is somewhat different in the infinite-measure category.

However, there are still some simple properties which do pass to quasi-factors, even in the infinite-measure category:

Recall that $(X,\mathcal{B},\mu,T)$ is \emph{rigid} if there exists  a strictly increasing sequence $(n_k)_{k \ge 0}$,
such that
for all $f \in L^2(\mu)$, $T^{n_k}f \to f$ in $L^2(\mu)$.
$(X,\mathcal{B},\mu,T)$ is \emph{rigid} iff
for all $A\in \mathcal{B}$ of finite measure, $\mu(T^{-n_k}A \bigtriangleup A) \to 0$.

\begin{proposition}\label{prop:rigid}
Suppose $(X,\mathcal{B},\mu,T)$ is rigid, then so is any quasi-factor.
\end{proposition}
\begin{proof}
Suppose $(n_k)_{k \in \mathbb{N}}$ is a rigid sequence for $T$, and let $(X^*,\mathcal{B}^*,\mu^*,T_*)$ .
It is sufficient for us to prove $T_*^{n_k}f \to f$  for a dense set of function  in $L^2(\mu^*)$.
Also, for bounded  $L^1$-functions, convergence in $L^1$ implies $L^2$-convergence.
For a measurable function $f:X \to \mathbb{R}$, Let $f^*:X^* \to \mathbb{R}$ be defined by
$$f^*(\gamma)=\int_{X}f(x)d\gamma(x).$$
Now,
$$\int_{X^*}| T^{n_k}_*f^*(\gamma)-f^*(\gamma)d\mu^*(\gamma)
 \le \int_{X^*}| T^{n_k}f-f|^*(\gamma)d\mu^*(\gamma) =$$
$$ = \int_X |T^{n_k}f(x)-f(x)|d\mu(x) \to 0$$
If $\int_X (T^{n_k}f-f)^2 d\mu \to 0$

Now the ring of functions generated by $\{ f^* :~ f \in L^1(\mu)\cap L^{\infty}(\mu)\}$
is dense in $L^2(\mu^*)$. Since $T^{n_k}_*h \to h$ for any function $h$ in this ring, our proof is complete.
\end{proof}

%
%

\section{\label{sec:zero_entropy}Quasi-factors of zero Krengel-entropy systems and zero entropy extensions}
Glasner and Weiss \cite{glasner_weiss95} proved that a probability preserving dynamical system
with zero entropy admits no quasi-factors with positive entropy.

It is natural to ask: Is it true that $0$ entropy for a \cmpts implies that any quasi-factor has zero entropy?
The answer, somewhat surprisingly is negative:

\begin{proposition}\label{prop:zero_krengel_pos_quasifactor}
 The exist \cmpts s with zero entropy in the sense of Krengel
(hence in the sense of Parry and Danilenko-Rudolph), which admit quasi-factors of positive entropy.
\end{proposition}
\begin{proof}
By proposition \ref{prop:factor_is_quasifactor} every non-singular factor $(\Omega,\mathcal{F},p,\tau)$,
which is a probability preserving transformation is isomorphic to a quasi-factor.
According to Ornstein and Weiss \cite{ornstein_weiss_infty_factor_84}, every non-singular transformation
is a non-singular factor of a tower over any other non-singular transformation.
Thus, we can begin with any ergodic zero-entropy probability preserving transformation $T_0$ (say a Kronecker system), construct a tower over $T_0$
that will admit any specific positive-entropy probability preserving transformation as a factor. The tower over $T_0$ will be a \cmpts with zero Krengel entropy,
which admits a quasi-factor with positive entropy.
\end{proof}

The result of \cite{glasner_weiss95} mentioned in the beginning of this section was generalized
by Glasner Thouvenot and Weiss \cite{glasner_thouvenot_weiss_2000}, to show that whenever
$\pi:(X,\mathcal{B},\mu,T) \to (Y,\mathcal{C},\nu,S)$  is a factor map of
probability preserving transformations such that the relative entropy $h(T \mid S)$ is zero,
and $(X^*,\mathcal{B}^*,\lambda,T_*)$ is a quasi-factor such that $\pi_*(\lambda) = \int_{Y^*}\delta_y d\nu(y)$, then $h_{\lambda}(T_* \mid \pi_*)=0$.
By taking $(Y,\mathcal{C},\nu,S)$ to be the trivial one-point system and $(X,\mathcal{B},\mu,T)$ to be a zero-entropy system,
this is a generalizes the previous result.

We prove the following theorem, which is a further generalization:

\begin{theorem}\label{thm:quasi-factor-relative-entropy}
Suppose $\pi:X \to Y$ is a zero entropy extension of conservative measure preserving
transformations $(X,\mathcal{B},\mu,T)$ and $(Y,\mathcal{C},\nu,S)$. Any quasi-factor $(X^*,\mathcal{B}^*,\xi,T_*)$ of $T$ is a
zero entropy extension of $\pi_*(\xi)$.
\end{theorem}

Theorem $7$ of \cite{glasner_thouvenot_weiss_2000} is a specific case of theorem \ref{thm:quasi-factor-relative-entropy}
where taking $X$ and $Y$ to be probability preserving systems,
and $\xi$ to be a quasi-factor such that $\pi_*(\xi)$ is distributed among point-mass probabilities on $Y$.

Our proof of theorem \ref{thm:quasi-factor-relative-entropy} is
analogous to Glasner and Weiss's proof of the corresponding theorem
$1.1$ in \cite{glasner_weiss95}, combined with ideas of Danilenko and Rudolph from \cite{danilenko_rudolph07}
 concerning relative entropy for infinite measure preserving transformations.

Recall that a factor map $\pi:X \to Y$ has relative completely
positive entropy if any for any intermediate factor $\pi_1:X \to Z$
such that $\pi_2:Z \to Y$ is also a factor and $\pi=\pi_2 \circ
\pi_1$, $Z$ has positive relative entropy. The \emph{relative Pinsker factor} $\mathcal{P}(T \mid S)$ is the smallest intermediate extension of $S$ ,
with respect to which $T$ has completely positive entropy. Equivalently: it is the largest factor of $T$ containing $S$, with respect to which
$T$ has zero relative entropy.

The following proposition appears in \cite{danilenko_rudolph07} as proposition $2.7$:

\begin{proposition}\label{prop:relative_entropy_join_dr}
If $(X_1,\mathcal{B}_1,\mu_1,T_1)$ and $(X_2,\mathcal{B}_2,\mu_2,T_2)$ are conservative systems with a common factor $(Y,\mathcal{B},\nu,S)$,
$h(T_1 \mid S) >0$ and $h(T_2 \mid S)=0$ then  $T_1$ and $T_2$ are relatively disjoint over $\mathcal{P}(T_1 \mid S)$, i.e. any joining of
$T_1$ and $T_2$ is relatively independent over $\mathcal{P}(T_1 \mid S)$.
\end{proposition}

\begin{lemma}\label{lem:rel_indenpedent_join2}
  Suppose $(Y,\mathcal{C},\nu,S)$ is a common factor of $(X_1,\mathcal{B}_1,\mu_1,T_1)$ and $(X_2,\mathcal{B}_2,\mu_2,T_2)$,
and suppose $T_1$ is a zero entropy extension of $S$.
Then the relatively independent joining of $T_1$ and $T_2$ over $S$ is a zero entropy extension of $T_2$.
\end{lemma}
\begin{proof}
  The relative entropy of $T_1 \times_{S} T_2$ over $T_2$ is the supremum over all finite partitions $\alpha$ of
  $h_\lambda(T_1\times T_2,\alpha \mid \mathcal{B}_2):=\lim_{n \to \infty}\frac{1}{n}H_\lambda( \alpha^n  \mid \mathcal{B}_2)$, where
  $$\alpha^n := \bigvee_{k=0}^{n-1} (T_1 \times T_2)^{-k}\alpha.$$
 By continuity of entropy (see \cite{danilenko_rudolph07} for a proof in the infinite-measure setting), it is sufficient to take $\alpha$ to be of the form
 $\alpha = \alpha_1 \vee \alpha_2$ where $\alpha_i$ is $\mathcal{B}_i$-measurable. For such $\alpha$'s,
  $$H_\lambda( \alpha^n \mid \mathcal{B}_2) = H_{\mu_1}(\alpha^n_1 \mid \mathcal{B}_2) + H_{\mu_2}( \alpha^n_2 \mid \mathcal{B}_2)$$
  Since $\alpha_2$ is $\mathcal{B}_2$-measurable, $H_{\mu_2}( \alpha^n_2 \mid \mathcal{B}_2)=0$.

  And so $h_\lambda(T_1\times T_2,\alpha) = h_{\mu_1}(T_1,\alpha_1 \mid \mathcal{B}_2) \ge h_{\mu_1}(T_1,\alpha_1 \mid \mathcal{C}) =0$.
\end{proof}

\textbf{Proof of theorem \ref{thm:quasi-factor-relative-entropy}:}
Consider system $(X\times X^*,\mathcal{B} \otimes \mathcal{B}^*, \tilde{\mu},T \times T_*)$, where $\tilde{\mu}$ is defined by:
$$\tilde{\mu}= \int_{X^*}\gamma \times \delta_{\gamma} d\mu^*(\gamma).$$

Similarly, consider $(Y\times Y^*,\mathcal{C} \otimes \mathcal{C}^*, \tilde{\nu},S \times S_*)$ with:
$$\tilde{\nu}= \int_{Y^*}\gamma \times \delta_{\gamma} d\nu^*(\gamma)$$

These systems are respectively, extensions of $T$ and $S$.

We have the following commuting diagram of homomorphism:

\xymatrix{ & X \times X^*,\tilde{\mu} \ar[dl]^{id \times \pi_*} \ar[dr]_{ \pi \times id} \ar[dd]_{\pi \times \pi_*}\\
X\times Y^*,\mu_1 \ar[dr]^{\pi \times id} & & Y \times X^*,\mu_2\ar[dl]_{id \times \pi_*} \\
& Y\times Y^*,\tilde{\nu} &}

 From this diagram it follows that $\mu_1$ is the relatively independent joining of $\mu$ and $\tilde{\nu}$ over $\nu$.
 Since $(X,\mathcal{B},\mu,T)$ is a zero entropy extension of $(Y,\mathcal{C},\nu,S)$, it follows from lemma \ref{lem:rel_indenpedent_join2}
that $(X \times Y^*,\mathcal{B} \otimes \mathcal{C}^*, \mu_1,T \times S_*)$ is a zero entropy
extension of $(Y\times Y^*,\mathcal{C} \otimes \mathcal{C}^*,\tilde{\nu},S \times S_*)$.

 For $A \in \mathcal{B}(Y)$ and $A^* \in \mathcal{B}^*$,
 $$\mu_2( A \times A^*  \mid \, (y,\pi^*(\gamma)))(y,\gamma) = 1_A(x)\mu^*(
 A^* \mid \, \pi^*(\gamma))(\gamma)$$

Assume for a moment that
 $(X^*,\mathcal{B}^*,\mu^*,T_*)$ is a positive entropy extension of $(Y^*,\mathcal{C}^*,\nu^*,S_*)$,
it follows that
 $(Y\times X^*,\mu_2)$ is a positive entropy extension of $(Y \times Y^*,
\tilde{\nu})$. In fact, in this case it is an infinite-entropy extension.

Thus, to prove the theorem it is sufficient to show that $\mu_2$ is
a zero entropy extension of $\tilde{\nu}$:

Suppose $(Y \times Z, \mathcal{C} \otimes \mathcal{Z},\mu_0,S \times S_0$) is the relative Pinsker factor of
$\pi_*:(Y \times X^*,\mu_2)\to( Y \times Y^*,\tilde{\nu})$, and
$p_1: Y\times X^* \to Y \times Z$, $p_2: Y\times Z \to Y \times Y^*$
are the corresponding factor maps.

It follows from proposition \ref{prop:relative_entropy_join_dr},
that under the conditional measure $\tilde{\mu}$ given $Y\times Z$, $X$ is independent of $X^*$:

$$\tilde{\mu}(\cdot  \mid ~ \pi(x)=y,\, \gamma)(x,\gamma)= \tilde{\mu}(\cdot  \mid ~ \pi(x)=y,\, p_1(y,\gamma))(x,\gamma)$$
But by definition of $\tilde{\mu}$,
$$\tilde{\mu}(\cdot \mid ~ \pi(x)=y, \gamma)(x , \gamma)= \gamma(\cdot \mid \pi(x)=y)(x)$$
and so the conditionals $\tilde{\mu}(\cdot \mid ~ \pi(x)=y,\, \gamma)(x,\gamma)$
 along with $\pi_*(\gamma)$ actually determine $\gamma$. In other words, the pullback of $\mathcal{C} \otimes \mathcal{Z}$
 is equal to $\mathcal{C} \otimes \mathcal{B}^*$.
 It follows that $p_1$ is an isomorphism, so $(Y \times X^*, \mu_2)$ is
its own Pinsker factor relative to $(Y\times Y^*, \tilde{\nu})$.
This completes the proof of the theorem.

 {\hfill $\Box$ \medskip\\}

\section{Positive relative entropy and quasi-factors}
Quasi-factors of probability preserving transformations with
positive entropy were considered in \cite{glanser_wiess_02}. The
main result is that each ergodic probability preserving system of
positive entropy admits every other ergodic probability preserving
system of positive entropy  as a quasi-factor. In this section we
 prove a ``relative'' generalization, and derive the following infinite-measure analog of this result:
\begin{theorem}\label{thm:positive_dan_rud_entropy_quasi_factors}
A conservative ergodic measure-preserving transformation with
positive Danilenko-Rudolph entropy admits any probability preserving
transformation as a quasi-factor.
\end{theorem}
The above theorem follows directly from the theorem below, by noting
that any system admits the trivial one-point system as a
quasi-factor:

\begin{theorem}\label{thm:positive_rel_entropy_quasifactor}
Given a positive entropy extension $T$ of a conservative ergodic
measure preserving transformation $S$,
   every positive entropy extension of any quasi-factor $\lambda_S$ of $S$ is relatively over isomorphic over $\lambda_S$ to a quasi-factor of $T$. In other words:

  Suppose the factor map $\pi:(X,\mathcal{B},\mu,T) \to (Y,\mathcal{C},\nu,S)$ has positive relative entropy and $(Y^*,\mathcal{C}^*,\lambda_S,S_*)$ is a quasi-factor of
  $(Y,\mathcal{C},\nu,S)$. Let $(\Omega,\mathcal{F},p,\tau)$ be a probability preserving transformation  such that
  $$\alpha: (\Omega,\mathcal{F},p,\tau) \to (Y^*,\mathcal{C}^*,\lambda_S,S_*)$$ is a factor map with positive relative entropy.
  Then there exist a quasi-factor $(X^*,\mathcal{B}^*,\lambda_T,T_*)$ of $(X,\mathcal{B},\mu,T)$ and an isomorphism
  $$\phi:(X^*,\mathcal{B}^*,\lambda_X,T_*) \to (\Omega,\mathcal{F},p,\tau),$$
  such that $\pi_* = \alpha \circ \phi$.

\end{theorem}

Danilenko and Rudolph asked if for every c.m.p.t the Krengel entropy is equal to the supremum of the relative entropies over all factors.
A positive answer to this question would imply along with theorem \ref{thm:positive_rel_entropy_quasifactor}
that a system with positive Krengel entropy admits every positive entropy p.p.t as a quasi-factor.

Let us state some results we need for the proof of theorem \ref{thm:positive_rel_entropy_quasifactor}:
\begin{lemma}\label{lem:quasifactor_of_factor}(Lemma $3.1$ of \cite{glanser_wiess_02})
If
$$\pi:(X,\mu) \rightarrow (Y,\nu)$$ is a homomorphism of m.p.t's,
then every quasi-factor $\lambda_0$ of $(Y,\nu)$ lifts to an isomorphic quasi-factor $\lambda$ of $(X,\mu)$
so that $\pi_*$ is an isomorphism of $\lambda$ and $\lambda_0$.
\end{lemma}
\begin{proof}
 Let
$$ \mu = \int_Y \mu_y d\nu(y)$$ be the disintegration of $\mu$ over $\nu$. The map $L:Y^* \to X^*$ defined by
$$L(\gamma)= \int_Y \mu_y d\gamma(y)$$
sends any quasi-factor $\lambda_0$ of $(Y,\nu)$ isomorphically  to the quasi-factor $\lambda := L_* \lambda_0$ of $(X,\mu)$,
and $\pi_* \lambda = \lambda_0$.
\end{proof}

In particular,
If
$$(X_2,\mu_2) \rightarrow_{\pi_1} (X_1,\mu_1) \rightarrow_{\pi_0} (X_0,\mu_0)$$ are homomorphism of m.p.t's and $\lambda_1$ is a quasi-factor of $(x_1,\mu_1)$ then it can be lifted to a quasi-factor

\begin{proposition}\label{prop:relative_Sinai}(Infinite-measure Sinai theorem)

  If $\pi:(X,\mu) \rightarrow (Y,\nu)$ is a homomorphism of ergodic \cmpts s  then there exists a m.p.t $(X_0,\mu_0)$ such that
  $(X_0,\mu_0) \to (Y,\nu)$ is relatively Bernoulli and
  $(X,\mu) \to (X_0,\mu_0)$ is relatively zero entropy.
\end{proposition}

\begin{proof}
   When the systems are probability preserving this proposition follows immediately from Sinai's theorem.
   To obtain the corresponding result for infinite measure preserving systems, choose some set of finite measure $A \in \mathcal{C}$,
    and apply Sinai's theorem to the induced transformation on $A$.
\end{proof}

Similarly, by inducing onto a set of finite measure, we obtain the
following variation on the Smorodinsky-Thouvenot theorem of
\cite{smord_thouv_1979}:
\begin{proposition}\label{prop:inf_smorodinsky} (Infinite-measure Smorodinsky-Thouvenot theorem)
Suppose $\pi:(X,\mathcal{B},\mu,T) \to (Y,\mathcal{C},\nu,S)$ is a factor map of positive entropy, then there exist $3$ factors of $ (X,\mathcal{B},\mu,T)$ such that the corresponding $\sigma$-algebras $\mathcal{B}_{i}$ $i=1,3,$ satisfy:
\begin{itemize}
  \item{ $\pi^{-1}\mathcal{C} \subset \mathcal{B}_i$ for $i=1,2,3$.}
  \item{ $\pi_i: (X,\mathcal{B}_i,\mu,T) \to (Y,\mathcal{C},\nu,S)$ is relatively Bernoulli.}
  \item{ These factors span: $\mathcal{B}= \bigvee_{i=1}^3 \mathcal{B}_i$.}
\end{itemize}
\end{proposition}
\begin{proof}
 The Smorodinsky-Thouvenot theorem of \cite{smord_thouv_1979}
 states that any probability preserving transformation is spanned by three Bernoulli factors.
 Namely, when $\mu(X) < \infty$, there exist $\sigma$-algebras $\mathcal{D}_i$, $i=1,2,3$ corresponding to
 Bernoulli factors so that $\mathcal{B}=\bigvee_{i=1}^3\mathcal{D}_i$.
 By Thouvenot's relative isomorphism theorem, $\mathcal{D}_i \vee \pi^{-1}\mathcal{C}$
 is a relatively Bernoulli extension of $S$. This proves the proposition for probability preserving systems.

 For infinite-measure conservative systems,
 the results is obtained by inducing onto a set $A \in \mathcal{C}$ of finite measure,
 and applying proposition \ref{prop:rel_isomorphism_induced} and lemma \ref{lem:rel_bernoulli_induce}.
\end{proof}

\noindent \textbf{Proof of theorem \ref{thm:positive_rel_entropy_quasifactor}:}
By Lemma \ref{lem:quasifactor_of_factor} and Proposition \ref{prop:relative_Sinai},
it is sufficient to prove the theorem in case $(X,\mathcal{B},\mu,T))$ is a relatively Bernoulli  extension of $(Y,\mathcal{C},\nu,S)$.

First we show that if $\pi: X \to Y$ is relatively Bernoulli $h(T
\mid S)= h_0$, and $\lambda_S$ is a quasi-factor of $S$, then there
is a quasi-factor $\lambda_T$ with $\pi_* \lambda_T =\lambda_S$ and
such that this is a relatively Bernoulli extension with relative
entropy $h$: There exist a set $A \in \mathcal{C}$ of finite
measure, and a $\mathcal{B}$-measurable partition
$\alpha=\{A_1,\ldots,A_N\}$ of $A$ with $H_{\mu}(\alpha)=h_0$
($\mu(A_i)=m_i$, $-\sum_{i=1}^{N}m_i \log(\frac{m_i}{\nu(A)})=h_0$)

 so that the iterates of $\alpha$ under $T_{A}$ are jointly independent, and also independent of $\mathcal{C}\cap A$,
and so that $\alpha$ is a relative generator with respect to the factor $S$:
$$\mathcal{B} \cap A= \bigvee_{n \in \mathbb{Z}}T^n_{A}\alpha \vee \mathcal{C} \mod \mu$$
To define $\lambda_T$, it is sufficient to describe the distribution of the measure of sets of the form $\bigcap_{k=1}^{n}T^{t_k}A_{j_k}$
with $A_{j_k} \in \alpha$, conditioned on the measure of sets in $\mathcal{C}$.
 We do so by requiring that all the mass in  $T^{-n}A$ will be concentrated in some $T^{-n}A_i \in T^{-n}\alpha$,
 with probabilities proportional to the measure of $A_i$, and that these events will be jointly independent. In other words:
$$\xi_T\left( \{\gamma : ~ \gamma(\bigcap_{k=1}^{n}T^{t_k}A_{j_k})= m\} \mid \{\gamma(\bigcap_{k=1}^{n}T^{t_k}A)= m \}\right)
 = \nu(A)^{-k}\prod_{k=1}^n m_{j_k}$$

Next, show that any $k$-fold self-joining over $\xi_Y$ of quasi-factors of the form $\xi_X$ is isomorphic to a quasi-factor $\lambda$ with
$\pi_*(\lambda)=\xi_Y$:

Suppose $((X^*)^k, \mathcal{B}^{\otimes k},\lambda, T_*^{\times k})$ is a $k$-fold self-joining of $\xi_X$, such that $\pi \circ p_i =\pi \circ p_j$,
where $p_i:(X^*)^k \to X^*$ is the projection onto the $i$'th coordinate.

Choose $k$ distinct positive numbers $\alpha_1,\ldots,\alpha_k$ with $\sum_{i=1}^k \alpha_i = 1$,
such that every subset-sum uniquely determines the subset. Namely, if $S, S' \subset \{1,\ldots,k\}$ are such that
$\sum_{i \in S} \alpha_i = \sum_{i \in S'}\alpha_i$,
  then $S= S'$.

The map $\phi:(X^*)^k \to X^*$ given by $\phi(\gamma_1,\ldots,\gamma_k)= \sum_{i=1}^k \alpha_i \gamma_i$ is $1-1$ almost-surely,
since the $\gamma_i$'s can be recovered uniquely:
$\gamma(A_j) = \sum_{i \in S}\alpha_i\gamma(A)$ iff $\gamma_i(A_j)= \gamma(A)$ for $i \in S$ and $\gamma_i(A_j)=0$ for $i \not\in S$.
It remains to check that $\phi(\lambda)$ is indeed a quasi-factor:
$$\int \gamma d\phi(\lambda)(\gamma)=$$
$$\int \sum_{i=1}^k\gamma_i d\lambda(\gamma_1,\ldots,\gamma_k) =
\sum_{i=1}^k \alpha_i\int \gamma_i d\xi_X(\gamma_i) = \mu$$

Now let $(\Omega,\mathcal{F},p,\tau)$ be any probability preserving
system which is an extension of a quasi-factor
$(Y^*,\mathcal{C}^*,\xi_Y,S_*)$ via $\rho: \Omega \to Y^*$. By the
infinite-measure Smorodinsky-Thouvenot theorem (proposition
\ref{prop:inf_smorodinsky}), There exists $\tau$-invariant
$\sigma$-algebras $\mathcal{F}_i \subset \mathcal{F}$ and
$\rho^{-1}\mathcal{C}^* \subset \mathcal{F}$ such that $\bigvee_i
\mathcal{F}_i = \mathcal{F}$ and such that Each factor
$\mathcal{F}_i$ is a relatively-Bernoulli extension of $\xi_Y$ with
relative entropy at most $h_0$. By our previous arguments,
 $(\Omega,\mathcal{F},p,\tau)$ is relatively-isomorphic (over $\xi_Y$) to a quasi-factor of $T$.
{\hfill $\Box$ \medskip\\}

Taking $\xi_Y$ in the above theorem to be the trivial quasi-factor (a $1$-point system),
 the following corollary regarding quasi-factors and
 Danilenko-Rudolph's entropy for infinite measure preserving transformations:
\begin{corollary}
  Any conservative system $(X,\mathcal{B},\mu,T)$ with $\hrel(T) >0$ admits any probability preserving system as a quasi-factor.
\end{corollary}

\section{Poisson suspensions}


The associated Poisson measure $\mu^* \in
\mathcal{P}(X^*,\mathcal{B}^*)$ is a probability measure which is
uniquely defined (up to isomorphism) by requiring that the random
measure of each set $A\in\mathcal{B}$ is distributed Poisson with
expectancy $\mu(A)$, and that the measures of disjoint sets are
independent:
\begin{equation}\label{eq:poisson_def}
\mu^*\left(\cap_{i=1}^n \{\gamma(A_i)=k_i\right\})=\prod_{i=1}^n
e^{-\mu(A_i)}\frac{\mu(A_i)^{k_i}}{k_i!}
\end{equation}
whenever $A_1,\ldots A_n \in \mathcal{B}$ are pairwise-disjoint.

Equation \eqref{eq:poisson_def} clearly implies that the barycenter
of $\mu^*$ is $\mu$. Also, $T_*$ preserves the Poisson-measure,
since equation \eqref{eq:poisson_def} and the independence of the
measures of disjoint sets hold for $\mu\circ T_*^{-1}$. Thus
$(X^*,\mathcal{B}^*,\mu^*,T_*)$ is a probability preserving system,
which by definition is the \emph{Poisson Suspension} of
$(X,\mathcal{B},\mu,T)$. By the above formula, for every $A \in
\mathcal{B}$ with $\mu(A)<\infty$,
$$\int_{X^*} \gamma(A) d \mu^*(\gamma) = \sum_{k=0}^{\infty}e^{-\mu(A)}\frac{\mu(A)^k}{k!}= \mu(A).$$
Thus, the Poisson-suspension $(X^*,\mathcal{B}^*,\mu^*,T_*)$ is a quasi-factor of $(X,\mathcal{B},\mu,T)$.
Poisson-suspensions of measure preserving systems have been studied by Roy \cite{roy_poisson}.
Various ergodic properties of the underlying infinite measure  preserving systems are reflected by ergodic properties of the Poisson suspension.
In some cases, this is  a mere consequence of the fact that the Poisson-suspension is a quasi-factor.

Applying theorem \ref{thm:quasi-factor-relative-entropy} to Poisson-suspensions immediately gives the
following:

\begin{corollary}\label{cor:poisson-zero-entropy-ext}
If $\pi:X \to Y$ is a factor map of the conservative measure
preserving system $(X,\mathcal{B},\mu,T)$ onto
$(Y,\mathcal{C},\nu,T)$, and $(X,\mathcal{B},\mu,T)$ zero
Krengel-entropy relative to $\pi$, then $(X,\mathcal{B},\mu,T)$  and
$(Y,\mathcal{C},\nu,T)$ have the same Poisson-entropy.
\end{corollary}

In \cite{jmrr_poisson_07} it was proved that the relative entropy of the
corresponding Poisson suspensions is equal to the relative entropy of the underlying transformations,
which gives a direct proof to corollary \ref{cor:poisson-zero-entropy-ext}.
However, for any $c>0$, the quasi-factor $\xi$ given by:
$$\xi \left(\gamma(A)=\frac{k}{c}\right)=e^{-c \cdot \mu(A)}\frac{c \cdot \mu(A)^k}{k!}$$
is isomorphic to the Poisson-suspension $(X^*,\mathcal{B}^*,(c\cdot\mu)^*,T_*)$  of $(X,\mathcal{B},c\cdot \mu,T)$
and by linearity of Poisson-entropy (also proved in \cite{jmrr_poisson_07}), the relative entropy $h_{\xi}(T_* \mid S_*)$ is $c$
times the relative entropy $h_\mu(T \mid S)$.

Another example of an application to Poisson-suspensions is proposition \ref{prop:rigid}, which implies that the Poisson suspension
of a rigid transformation is rigid - a result obtained in theorem $4.3$ of \cite{roy_poisson}.

On the other hand, many properties exhibited by Poisson suspensions are not shared by all quasi-factors.

In this context, we mention a question raised in \cite{jmrr_poisson_07}: Is it true that the Poisson suspension of any \cmpts
with zero Krengel entropy has zero entropy? As we have seen this is not the case with general quasi-factors.

%
%

\bibliographystyle{abbrv}
\bibliography{quasifactors}
\end{document}